\documentclass{article}
\usepackage{ijcai26}

\usepackage{times}
\usepackage{soul}
\usepackage{url}
\usepackage[utf8]{inputenc}
\usepackage[small]{caption}
\usepackage{graphicx}
\usepackage{amsmath,amssymb,amsfonts}
\usepackage{amsthm}
\usepackage{booktabs}
\usepackage[ruled,vlined]{algorithm2e}
\usepackage[switch]{lineno}
\usepackage{bm}
\usepackage{xcolor}
\usepackage[
  colorlinks=true,
  linkcolor=blue,
  citecolor=blue,
  urlcolor=blue
]{hyperref}
\usepackage[nocomma]{optidef}

\usepackage{multirow}
\usepackage{array}
\usepackage{tikz}
\usepackage{pgfplots}
\pgfplotsset{compat=1.18}

\title{Column Generation for the Micro-Transit Zoning Problem}

\author{
Hins Hu$^1$
\and
Rishav Sen$^2$\and
Jose Paolo Talusan$^2$\and
Abhishek Dubey$^2$\and
Aron Laszka$^3$\and
Samitha Samaranayake$^1$
\affiliations
$^1$Cornell University\\
$^2$Vanderbilt University\\
$^3$Pennsylvania State University\\
\emails
zh223@cornell.edu
}
\begin{document}

\maketitle

\begin{abstract}
Along with the rapid development of new urban mobility options like ride-sharing over the past decade, on-demand micro-transit services stand out as a middle ground, bridging the gap between fixed-line mass transit and single-request ride-hailing, balancing ridership maximization and travel time minimization. Micro-transit adoption can have significant social impact. It improves urban sustainability, through lower energy consumption and reduced emissions, while enhancing equitable mobility access for disadvantaged communities, thanks to its lower vehicle miles per passenger, flexible schedules, and affordable pricing. However, effective operation of micro-transit services requires planning geo-fenced zones in advance, which involves solving a challenging combinatorial optimization problem. Existing approaches enumerate candidate zones first and selects a fixed number of optimal zones in the second step. In this paper, we generalize the Micro-Transit Zoning Problem (MZP) to allow a global budget rather than imposing a size limit for candidate zones. We also design a Column Generation (CG) framework to solve the problem and several pricing heuristics to accelerate computation. Extensive numerical experiments across major U.S. cities demonstrate that our approach produces higher-quality solutions more efficiently and scales better in the generalized setting.
\end{abstract}

\section{Introduction}
\label{sec:intro}
Public transit has long been fundamental to urban mobility systems, offering high energy efficiency \cite{hodges2010public}, reduced emissions, substantial capacity, and social equity benefits \cite{miller2016public}. However, its lengthy planning horizons, significant infrastructure investments, and budgetary constraints have resulted in limited coverage across most U.S. cities, leaving disadvantaged communities with poor accessibility to essential services and extended travel times. For example, according to the \href{https://www.census.gov/topics/employment/commuting/guidance/acs-1yr.html}{American Commuter Survey}, only 3.7\% of commuters in the U.S. use public transit for work commute in 2024.

On the other hand, on-demand ride-sharing has gained popularity as an alternative over the past decade, providing responsive door-to-door services \cite{kooti2017analyzing}. Yet, it presents critical drawbacks: prohibitive costs for low-income commuters, higher per-passenger emissions, and insufficient pooling that increases deadheading, which exacerbates congestion and undermines sustainability.

\textit{Micro-transit} has recently emerged as a middle ground between these two modes \cite{bardaka2024empathy}. This relatively high-capacity service mixes pre-booked and on-demand requests, and combines flexible routing and adaptive scheduling, filling coverage gaps left by fixed-route transit. It serves both short neighborhood trips and provides first-mile/last-mile connections to major hubs. With affordable pricing, micro-transit balances travel time, coverage, and cost across the urban mobility spectrum.

The first crucial technical challenge for city planners operating micro-transit is determining where to establish geo-fenced service areas given existing monetary budgets, specific requirements from local communities, and operational objectives (e.g., maximizing demand served). Micro-transit zoning design is a complex decision-making problem typically addressed by human experts without algorithmic processes or formal modeling. While various computational zoning approaches exist, many are adapted from different contexts, such as freight logistics \cite{chandra2020designing} and urban development.


To the best of our knowledge, \cite{hu2025optimal} proposed the first optimal micro-transit zoning method to maximize total demand coverage. This two-phase framework develops an iterative algorithm with a geographical heuristic to enumerate all valid candidate zones under a size limit in Phase 1, then solves a Weighted Maximum Covering Problem via integer programming in Phase 2 to identify $m$ optimal zones, where $m$ is a problem input.

However, this work has notable limitations: (1) the zoning model assumes prior knowledge of the number of zones to establish (i.e., $m$ is a fixed input), rather than allowing a global budget for operating all zones within an urban region's boundary, which is less realistic; (2) the enumeration algorithm in Phase 1 suffers from scalability issues in larger cities when demand aggregation becomes more granular (e.g., less than $1$ km$^2$).

In this paper, we propose a more efficient approach using Column Generation (CG) to address a generalized zoning design setting with a \textit{global budget}, which we call the Micro-Transit Zoning Problem (MZP). Rather than relying solely on geographical heuristics for candidate zone generation, we reformulate the problem into a CG framework consisting of a \textit{restricted master problem} and a \textit{pricing problem}, both modeled as Integer Linear Programs (ILPs). Under the CG framework, we leverage mathematical programming to handle exponentially many decision variables, aligning with the nature that the candidate set of valid zones is exponentially large. To further accelerate computation, we develop a pricing heuristic to replace the exact pricing ILP with negligible compromise in solution quality. The effectiveness of our framework is validated through extensive experiments across multiple large U.S. cities. Numerical results demonstrate the superiority of the proposed CG-based approach over the method by \cite{hu2025optimal}. Additional analysis and parameter tuning for our method are provided to better guide engineering practice.

\subsection{Stakeholder Collaboration and Societal Impact}
This research addresses a real-world societal challenge: the need for computationally tractable and equitable approaches to design micro-transit service zones for cities that face low transit accessibility and financial burdens to expand their fixed-line transit networks. The work emerges from an ongoing multi-year collaboration with the Chattanooga Area Regional Transportation Authority (CARTA), our primary stakeholder and domain partner, who serves a region where only 1.6\% of commuters currently use public transit for work commutes, significantly below the national average. 

CARTA actively participates in this project through three primary roles: (1) conducting community surveys and providing insights for algorithm design; (2) coordinating with other local stakeholders, such as social workers and mini-bus drivers,
to calibrate our model parameters; (3) leading a micro-transit \href{https://www.gocarta.org/}{pilot project} for real-world validation.

\subsection{Contributions}
Our contributions to the research community of AI for social good are three-fold:
\begin{enumerate}
    \item \textbf{Generalized Problem and Formulation:} We generalize the Micro-Transit Zoning Problem (MZP) to incorporate a global budget constraint on operational costs and reformulate to a Column Generation (CG) framework. 
    \item \textbf{Scalable Pricing Heuristic:} We develop an efficient pricing heuristic algorithm for scalability.
    \item \textbf{Real-World Validation:} We conduct extensive experiments across multiple large U.S. cities using real-world mobility demand data, demonstrating superior performance over state-of-the-art methods in both solution quality and computational efficiency.
\end{enumerate}

\section{Related Work}
In this section, we first review the existing computational methods for zoning in urban planning and transportation, including ones particularly for micro-transit services. Second, we briefly discuss the development of Column Generation (CG) and its applications on relevant problems in transportation, such as ride-sharing and vehicle routing. 

\subsection{Zoning in Urban Planning and Transportation}
The design of spatial zones plays a critical role in urban planning and infrastructure development, shaping land use decisions \cite{lens2022zoning} and influencing the trajectory of urban growth \cite{domingo2021effect}. Inadequately designed or obsolete zoning can intensify spatial inequalities in cities \cite{dingil2025fostering}. In the domain of general transportation planning, a great number of computational methods have emerged to support automated or optimal zoning. For example, approaches have been developed using K-means clustering \cite{chandra2020designing} or multi-objective genetic algorithms \cite{chandra2021multi} to refine zone structures by grouping geographically similar spatial units.

In the specific context of micro-transit zoning, computational methods remain limited. For example, \cite{liu2021mobility} proposed a framework integrating fixed-route transit with zone-based demand-responsive services, where the zone size is optimized jointly with transit network design. However, their approach assumes uniform zones adjacent to each other in a continuous region. \cite{bonner2023achieving} developed a mixed-integer programming approach for equitable micro-transit zone generation, prioritizing fair service distribution across populations and geographic areas. Most relevant to our work, \cite{hu2025optimal} formulated the micro-transit zoning problem to maximize demand coverage and proposed a two-phase solution framework involving candidate zone generation and optimal selection. As stated in Section \ref{sec:intro}, our work extends their problem setting and provides a more scalable solution approach based on Column Generation.

\subsection{Column Generation} \label{sec:review_cg}
Column Generation (CG) emerged in the 1950s as a technique for solving large-scale Linear Programs (LP), where explicitly enumerating all decision variables is computationally prohibitive. The seminal work of \cite{dantzig1960decomposition} established its theoretical foundations. Early successes in the Cutting Stock Problems \cite{gilmore1961linear} demonstrated its practical value. The integration of CG with branch-and-bound led to the development of branch-and-price algorithms \cite{barnhart1998branch}, enabling the solution of large-scale Integer Programs (IP). As of now, CG had become a cornerstone methodology in artificial intelligence, particularly for problems with decomposable structures. Further details can be found in \cite{desaulniers2006column}.

CG has been widely applied across various fields, including vehicle routing, scheduling, resource allocation, ride-sharing, and more. \cite{desrochers1992new} solved the Vehicle Routing Problem with Time Windows using a set-partitioning formulation where each column represents a feasible route, and the pricing problem is formulated as a resource-constrained shortest path problem. \cite{vance1997airline} applied CG to the Airline Crew Pairing Problem, generating crew parings (i.e., columns) that cover flight segments while satisfying labor regulations - a method now widely adopted in the airline industry. 
\cite{van1999parallel} employed CG for the Parallel Machine Scheduling Problem, where columns represent feasible job sequence for individual machines. \cite{riley2019column} developed a real-time share-a-ride system using column generation to handle up to 30,000 ride requests per hour in New York City. Our paper contributes to this repertoire by applying CG to the Micro-Transit Zoning Problem.

\section{Problem Statement}
\label{sec:problem}
We define the \textbf{Micro-Transit Zoning Problem (MZP)} as a combinatorial optimization problem with the following inputs:
\begin{enumerate}
    \item A set of uniform, non-overlapping cells $V$ in an urban region, with each cell representing a small neighborhood that aggregates all travel demand from and to the area. Non-adjacent cells are allowed if no travel demand arises in the area in between. \label{input:1}
    \item A road network $G$ in the same region where $V$ lies on.
    \item A demand table $\bm{d}: V^2 \rightarrow \mathbb{R}_+$, where each entry $d(i, j)$ represents a nominal travel demand from cell $i$ to $j$.
    \item A micro-transit zone $S \in 2^V$ is defined as a subset of cells that are \textit{adjacent} to each other in the \textit{road network space}, not the 2-D geographical space.
    \item A distance matrix $\bm{c}: V^2 \rightarrow \mathbb{R}_+$, where each entry $c(i, j)$ denotes the shortest-path distance in the network $G$ from the centroid of cell $i$ to the centroid of cell $j$. 
    \item A linear cost function $f(S) = \alpha \, D_S^2 + \beta$ that maps the \textit{squared diameter} $D_S^2$ of a zone $S$ to the cost of operating that particular zone. The diameter is defined as distance between the most distant two cells within a zone, namely $D_S \triangleq \max_{i, j \in S} c(i, j)$. \label{input:6}
    \item A global budget $B \in \mathbb{R}_+$ for total cost of operating all micro-transit zones. \label{input:7}
\end{enumerate}
The objective is to select a collection of micro-transit zones to operate under budget $B$, so as to maximize the total intra-zone demand shown in Equation (\ref{eq:obj}).
\begin{equation} \label{eq:obj}
\sum\limits_{i,j \in V} d(i, j) \; \mathbb{I} (i, j \text{ in the same zone})  
\end{equation}

The problem studied here generalizes that of \cite{hu2025optimal}. Specifically, we model the monetary constraint of transit agencies as a global budget for all operated zones collectively, rather than imposing a uniform size limit on each individual zone, resulting in a more realistic setting.

For input \ref{input:1}, non-uniform cells can be designed to match demand patterns without affecting our approach. Input \ref{input:6} reflects our partner transit agency's insight that operational cost is proportional to area size. Beyond the global budget (input \ref{input:7}), single-zone budget can also be incorporated if needed. 

\section{Methodology}
In this section, we briefly introduce Column Generation (CG) and its key terminology, then describe each component as applied to our micro-transit zoning problem.

\subsection{Prerequisite for Column Generation} \label{sec:cg_prerequisite}
CG is a powerful decomposition technique particularly effective for problems that can be formulated as \textit{set partitioning} or \textit{set covering} models, where each decision variable represents a partial solution with special structure - referred to as a \textit{column}.

Due to the exponentially growing nature of feasible columns with the problem size, CG works on a \textit{Restricted Master Problem (RMP)}, which has the same structure as the original problem but includes only a subset of columns. Starting with an initial small subset, CG iteratively identifies and adds promising columns to the RMP.

\begin{figure}[htbp!]
    \centering
    \includegraphics[width=\linewidth]{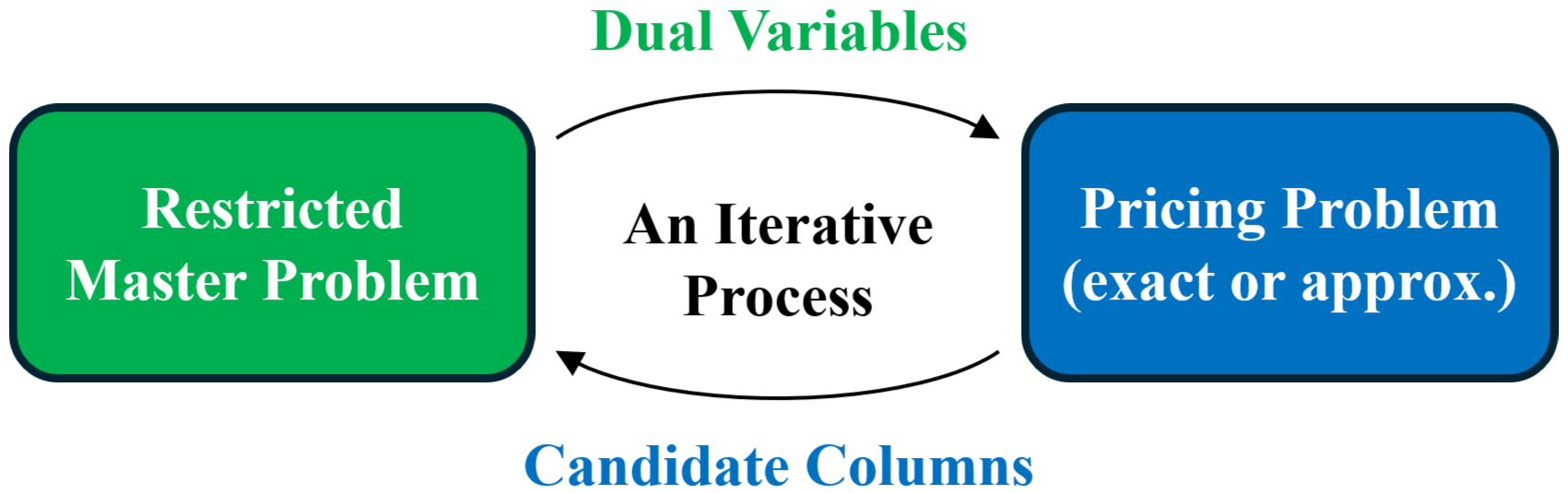}
    \caption{A Sketch of Column Generation Process}
    \label{fig:cg-scheme}
\end{figure}

As shown in Figure \ref{fig:cg-scheme}, the process alternates between two steps. First, the RMP is solved to obtain an optimal solution and associated \textit{dual variables}. Second, these duals are passed as input to the \textit{Pricing Problem} - an auxiliary optimization problem that searches for columns with negative \textit{reduced cost} \footnote{If the RMP is a maximization, the pricing aims for columns with positive reduced cost, which is exact our case.} Any column with negative reduced cost indicates a potential increase on the RMP objective and is added to the column set.

CG terminates when the pricing problem no longer identifies improving columns, guaranteeing optimality of the linear relaxation to the unrestricted master problem. To obtain integer solutions, column generation is typically embedded within a branch-and-bound framework, known as \textit{branch-and-price} \cite{barnhart1998branch}, where column generation solves the LP relaxation at each node of the search tree.

\subsection{The Restricted Master Problem}
\label{sec:master}
As discussed in Section \ref{sec:cg_prerequisite}, applying CG requires a set-based formulation. In the micro-transit zoning problem, each column naturally corresponds to a zone as it is a collection of cells, which is analogous to vehicle routes serving as columns in the Vehicle Routing Problems \cite{desrochers1992new}. 

In addition to the notations introduced in Section \ref{sec:problem}, we further define the following for the formulation of the CG framework. Let $\mathcal{S}$ be the set of candidate zones. Let $x_S \in \{0, 1\}$ be a binary decision variable indicating whether to select zone $S$. Let $w_{ij} \in \{0, 1\}$ be an auxiliary decision variable indicating whether two distinct cell $i$ and $j$ are in the same zone, which is used to avoid double-counting the demand whose origin and destination are both in the overlapping area of two distinct zones. Let $z_S = [z_S^1, \dots, z_S^i, \dots, z_S^{|V|}]^\top$ be the cell-selection vector of zone $S$, where $z_S^i$ is a binary coefficient indicating whether cell $i$ is included in zone $S$. Then, we formulate the RMP as in ILP (\ref{ilp:master}). 
\begin{subequations}
\label{ilp:master}
\begin{alignat}{2}
\max_{\bm{x},\bm{w}} \quad 
& \sum_{i,j \in V} d_{ij}\, w_{ij} \\[2mm]
\text{s.t.} \quad
& \sum_{S \in \mathcal{S}} (\alpha D_S^2 + \beta)\, x_S \le B 
& \quad (\lambda) \label{c:master_budget_dual} \\[1mm]
& w_{ij} - \sum_{S \in \mathcal{S}} z_S^i z_S^j\, x_S \leq 0
& \quad (\pi_{ij}) \quad \forall i,j \in V \label{c:master_demand_coverage} \\[1mm]
& x_S \in \{0,1\} 
& \quad \forall S \in \mathcal{S} \label{c:master_int_1}\\[1mm]
& w_{ij} \in \{0,1\} 
& \quad \forall i,j \in V \label{c:master_int_2}
\end{alignat}
\end{subequations}

Constraint (\ref{c:master_budget_dual}) requires that the cumulative operational cost of all selected zones be capped by the global budget. Constraint (\ref{c:master_demand_coverage}) links $\bm{x}$ and $\bm{w}$, forcing $w_{ij}$ to be $0$ if cell $i$ and $j$ are not assigned simultaneously to any zone $S$ (with at least one of $z_S^i$, $z_S^j$, and $x_S$ being $0$). In the maximization setting, it also raises $w_{ij}$ to $1$ if cell $i$ and $j$ are assigned to at least one zone (with $z_S^i = 1$, $z_S^j = 1$, and $x_S = 1$ for at least one $S \in L$). Constraints (\ref{c:master_int_1}) and (\ref{c:master_int_2}) enforce the integrality of decision variables. For better readability, the dual variables $\lambda$ and $\pi_{ij}$ are marked next to their corresponding constraints.

Despite its theoretical soundness, the above formulation often exhibits degeneracy in real-world instances, which can cause the CG process to stall or cycle \cite{du1999stabilized}. As an engineering practice, we perturb Constraint (\ref{c:master_demand_coverage}) by adding small random noise terms $\epsilon_{ij}$ to the right-hand side, which breaks degeneracy at the cost of non-monotonic convergence in the master objective value.
\begin{equation} \label{c:random_noise}
w_{ij} - \sum_{S \in \mathcal{S}} z_S^i z_S^j\, x_S \leq \epsilon_{ij} \quad \forall \; i,j \in V 
\end{equation}

According to \cite{bertsimas1997introduction}, the reduced cost for zone $S$ is can be computed by Equation (\ref{eq:reduced_cost}), which will be used in the pricing problem for determining whether $S$ can be added to the candidate column set.
\begin{equation}
\label{eq:reduced_cost}
\bar{c}_S = \sum_{i,j\in V} \pi_{ij} \, z_S^i \, z_S^j - \lambda \, (\alpha \, D_S + \beta)
\end{equation}

ILP (\ref{ilp:master}) shares the structure of \textsc{ZoningILP} in \cite{hu2025optimal} except for adopting a global budget. Notably, the candidate zone set $\mathcal{S}$ in our CG framework is orders of magnitude smaller, as dual variables guide the pricing problem to generate only ``premuim'' candidate zones.

\subsection{The Exact Pricing Problem}
In the pricing problem, the core idea is to find a zone to maximize the reduced cost, and thus the cell selection vector $z_S$ become part of the decision and the index $S$ for zones can be discarded. Formally, we define a binary decision variable $z_i$ for every cell $i \in V$ to indicate whether $i$ is selected in the zone.

Similarly, we define a positive continuous decision variable $D^2 \in \mathbb{R}_+$, with the index $S$ discarded as well, to represent the actual ``squared diameter'' of the selected zone. This quantity is an input to the RMP because the set of candidate zones are given. However, in the pricing problem, it varies in the optimization process guided by the dual variables.

Then, an exact Integer Quadratic Programming (IQP) formulation for pricing is shown in IQP (\ref{qp:pricing}). From Constraint (\ref{c:QP_pricing_diameter}), as long as both cell $i$ and $j$ are selected simultaneously in a zone, its diameter is lowered bounded by the two-way shortest path distance between them. The max operator is used to account for the asymmetry in city road networks. The objective function is adopted from the reduced cost in Equation (\ref{eq:reduced_cost}) by converting $\bm{z}$ and $D$ to decision variables and dropping the constant term $\lambda \, \beta$.
\begin{subequations}
\label{qp:pricing}
\begin{alignat}{2}
\max_{\bm{z},\, D^2} \quad
& \sum_{i,j \in V} \pi_{ij}\, z_i z_j - \lambda\, \alpha\, D^2 \\[2mm]
\text{s.t.} \quad
& D^2 \geq \max\{c^2(i,j), c^2(j, i)\}\, z_i z_j
& \quad \forall i,j \in V \label{c:QP_pricing_diameter} \\[1mm]
& z_i \in \{0,1\}
& \quad \forall i \in V
\end{alignat}
\end{subequations}

Since a generic IQP can be computationally challenging for modern mixed-integer programming solvers, we further linearize the bilinear term ($z_i \, z_j$) to obtain an ILP pricing problem as an alternative. Let $y_{ij}$ be a binary decision variable for each pair of cells $(i,j)$, indicating whether both cells are selected. This definition coincides with that of $\bm{w}$ in the master problem (Section \ref{sec:master}), although the variables serve different roles. Then, the linearized formulation is shown in ILP (\ref{ilp:pricing}).
\begin{subequations}
\label{ilp:pricing}
\begin{alignat}{2}
\max_{\bm{y}, \bm{z}, D^2} \quad
& \sum_{i, j \in V} \pi_{ij} \, y_{ij} - \lambda \, \alpha \, D^2 \\[2mm]
\text{s.t.} \quad
& y_{ij} \leq z_i
& \quad \forall i, j \in V \label{c:linearize_1} \\[1mm]
& y_{ij} \leq z_j
& \quad \forall i, j \in V \label{c:linearize_2}\\[1mm]
& y_{ij} \geq z_i + z_j - 1
& \quad \forall i, j \in V \label{c:linearize_3}\\[1mm]
& D^2 \geq \max\{c^2(i, j), c^2
(j, i)\} \, y_{ij}
& \quad \forall i, j \in V \\[1mm]
& y_{ij} \in \{0,1\}
& \quad \forall i, j \in V \\[1mm]
& z_{i} \in \{0,1\}
& \quad \forall i \in V
\end{alignat}
\end{subequations}

\subsection{Pricing Heuristic}
\label{sec:pricing_heuristic}

While the exact pricing formulations in IQP (\ref{qp:pricing}) and ILP (\ref{ilp:pricing}) guarantee optimality, they can still be computationally intractable for large-scale instances, particularly when the pricing problem must be solved repeatedly throughout the CG process. To ensure the effectiveness of the framework, we develop a \textit{construction heuristic} that rapidly generates high-quality zones with positive reduced cost.

Motivated by the generation of school bus routes in \cite{bertsimas2019optimizing}, our heuristic employs a greedy node-addition strategy: starting  from a randomly selected pair of cells (which should satisfy the single-zone budget constraint if enforced), it iteratively adds cells that maximize the marginal gain in reduced cost until no beneficial addition remain. The marginal gain of adding cell $k$ to the current zone $S$ can be efficiently computed by Equation (\ref{eq:marginal_gain}).
\begin{equation} \label{eq:marginal_gain}
\Delta(k) = \sum_{j \in S} (\pi_{kj} + \pi_{jk}) - \lambda \cdot \alpha \cdot (D^2_{\text{new}} - D^2_{\text{current}})
\end{equation}
where $D_{\text{new}}$ represents the updated squared diameter after including cell $k$, and $D_{\text{current}}$ is the current diameter of zone $S$. The first term captures the added demand coverage from new within-zone connections, while the second term accounts for the increased operational cost due to zone expansion.

To enhance robustness and exploration, the heuristic can be executed multiple times with different random initialization. Each run produces a local optimal zone from a particular seed pair of nodes. Upon termination, all generated zones with positive reduced costs are added to the RMP's candidate set $\mathcal{S}$. Algorithm \ref{alg:pricing_heuristic} presents the complete procedure.

\begin{algorithm}[htbp!]
\caption{Greedy Pricing with Random Seeds}
\label{alg:pricing_heuristic}
\KwIn{cell set $V$, distances $\bm{c}$, duals $\lambda$ and $\{\pi_{ij}\}$, parameters $\alpha$ and $\beta$, number of runs $R$, single-zone budget $B_0$}

\KwOut{list of zones $\mathcal{Z}$ with positive reduced cost}

$\mathcal{Z} \gets \emptyset$, 

\For{$r = 1$ \KwTo $R$}{
    
    \tcp{Random Initialization}
    \Repeat{$\alpha \cdot D^2 + \beta \leq B_0$}{
        Randomly select two distinct cells $i, j \in V$\\
        $D \gets \max\{\bm{c}[i,j], \bm{c}[j,i]\}$\\
    }
    $S \gets \{i, j\}$\\
    
    \tcp{Greedy Node Addition}
    \Repeat{convergence or timeout}{
        $k^* \gets \arg\max_{k \in V \setminus S} \Delta(k)$\\
        
        $\hat{D} \gets \max_{j \in S} \max\{\bm{c}[k^*,j], \bm{c}[j,k^*]\}$\\
        $D_{\text{new}} \gets \max\{D, \hat{D}\}$\\

        \If{$\Delta(k^*) > 0$ \textbf{and} $\alpha \cdot D^2_{\text{new}} + \beta \leq B_0$}{
            $S \gets S \cup \{k^*\}$\\
            $D \gets D_{\text{new}}$\\
        }
        \Else{
            \textbf{break}
        }
    }
    \tcp{Zone Selection}
    Compute $\bar{c}_S = \sum_{i,j \in S} \pi_{ij} - \lambda(\alpha D^2 + \beta)$\\
    \If{$\bar{c}_S > 0$}{
        $\mathcal{Z} \gets \mathcal{Z} \cup \{S\}$
    }
}
\Return $\mathcal{Z}$
\end{algorithm}

The heuristic offers several fundamental advantages. First, it executes in $O(|V|^2)$ time per random run, significantly faster than solving the exact pricing problem, either in the formulation of ILP or IQP. Second, it increases column diversity through random initialization, thus lowering the risk of the CG process being trapped in local optima. Third, the heuristic supports a time limit, functioning as an \textit{anytime} algorithm that enables flexible trade-offs between efficiency and solution quality.

\subsection{The Overall Framework}
We summarize the complete CG framework for the micro-transit zoning problem, as presented in Algorithm \ref{alg:cg_framework}. 
\begin{algorithm}[htbp!]
\caption{CG for Micro-Transit Zoning}
\label{alg:cg_framework}
\KwIn{cell set $V$, demand $\bm{d}$, distances $\bm{c}$,  parameters $\alpha$ and $\beta$, budget $B$}
\KwOut{Selected zones $\mathcal{S}^\star$}

\tcp{Initial Candidate Set}
$\mathcal{S} \gets \text{a random zone under single-zone budget $B_0$}$

\tcp{Column Generation}
\Repeat{$\mathcal{S}_{new} = \varnothing$ or \text{timeout}}{
$(\lambda, \{\pi_{ij}\}) \gets \textsc{RMP}(\mathcal{S})$

\If{\text{exact}}{ 
$\mathcal{S}_{new} \gets \textsc{ExactPricing}(\lambda, \; \{\pi_{ij}\}, \; \dots)$
}
\Else{
$\mathcal{S}_{new} \gets \textsc{PricingHeuristic}(\lambda, \; \{\pi_{ij}\}, \; \dots)$
}
$\mathcal{S} \gets \mathcal{S} \; \cup \; \mathcal{S}_{new}$
}
\tcp{Solve Restricted Master}
$\mathcal{S}^* \gets \textsc{RMP}(\mathcal{S})$

\Return $\mathcal{S}^\star$
\end{algorithm}

Our framework only applies CG in the LP relaxation to the RMP (i.e., the root node in the branch-n-bound tree), without embedding the procedure into a more sophisticated full branch-n-price scheme. While the latter would be theoretically required to guarantee optimality, its implementation would demand substantial re-engineering of lower-level, fine-tuned branching rules and cutting plane algorithms encapsulated in any modern Mixed Integer Programming (MIP) solver, making it an impractical option. 

Nevertheless, through the numerical experiments in Section \ref{sec:exp}, we show that our framework is sufficiently superior to the state-of-the-art baseline. 

\section{Numerical 
 Experiments} \label{sec:exp}
In this section, we present comprehensive experiments to evaluate the effectiveness and scalability of our CG framework on real-world micro-transit zoning instances. We consider five medium-to-large U.S. cities: Miami, Boston, Atlanta, Chattanooga, and Nashville.

\subsection{Data Preparation}
City boundaries and underlying road networks are publicly available. For Chattanooga, travel demand is derived from a proprietary dataset of real origin-destination (OD) trips provided by our partner transit agency, CARTA. For all other cities, demand is synthetically generated following the methodology in \cite{sen2025moveod}, which leverages two real-world sources: (1) departure time and travel time distributions from the American Community Survey (ACS), and (2) residence-to-workplace flows from the Longitudinal Employer-Household Dynamics (LODES). In pre-processing, we apply a geo-fencing filter to retain only trips with both origins and destinations inside the city boundaries and exclude very short trips with total travel distance less than 500 meters.
\begin{table}[htbp!]
\centering
\caption{H3 Hexagon and Road Network Statistics Across Multiple Cities}
\label{tab:graph_stats}
\resizebox{0.45\textwidth}{!}{%
\begin{tabular}{c cc cc}
\toprule
\multirow{2}{*}{\textbf{City}} 
& \multicolumn{2}{c}{\textbf{\# of H3 Hexagons}} 
& \multicolumn{2}{c}{\textbf{Road Network}} \\
\cmidrule(lr){2-3} \cmidrule(lr){4-5}
& \textbf{Res. 7} & \textbf{Res. 8} 
& \textbf{\# of nodes} & \textbf{\# of edges} \\
\midrule
Miami        & 31  & 220  & 9086  & 23558 \\
Boston       & 45  & 321  & 11408 & 26086 \\
Atlanta      & 73  & 502  & 12905 & 33628 \\
Chattanooga  & 80  & 551  & 9486  & 24358 \\
Nashville    & 257 & 1803 & 22470 & 55796 \\
\bottomrule
\end{tabular}%
}
\end{table}

To spatially aggregate travel demand into small, contiguous neighborhoods - each corresponding to a cell as defined in Section \ref{sec:problem} - we the open-source \href{https://github.com/uber/h3}{Hexagonal Hierarchical Spatial Index (H3)} developed by Uber. H3 divides geographic space into hexagonal grids at multiple discrete resolutions. We focus on resolution 7 (average cell length $\approx 1.22$ km) and resolution 8 ($\approx 0.46$ km), which provide a reasonable granularity for micro-transit zoning based on insights from CARTA. Trips are assigned to hexagons by mapping their origin and destination coordinates to the corresponding cells. As an example, the aggregated demand pattern in Chattanooga is visualized in Figure \ref{fig:out_demand_chatt}. Full details of visualization for other cities can be found in Appendix A.
\begin{figure}[htbp!]
    \centering
    \includegraphics[width=\linewidth]{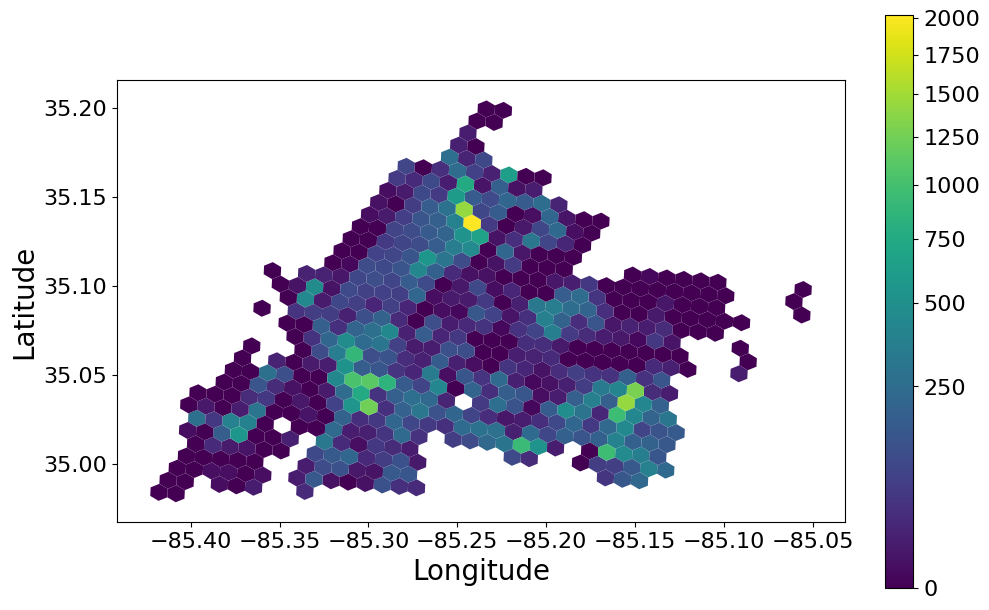}
    \caption{Aggregated number of trips in Chattanooga with H3 Resolution 8}
    \label{fig:out_demand_chatt}
\end{figure}

In summary, geographic statistics for all cities are reported in Table \ref{tab:graph_stats}, and the total demand can be found in Appendix C. 



\subsection{Experimental Settings}
\label{sec:exp_setup}
The experiments are conducted under the following configurations and parameter settings unless otherwise specified. We use pre-computed nominal driving times on the real-world road networks to construct the distance matrix $\bm{c}$. To make operational costs comparable across cities of different sizes, distances are normalized to the range $[0,1]$.

The cost parameters are set uniformly across all cities: global budget $B = 8$, cost scale factor $\alpha = 5$, and fixed cost per zone $\beta = 1$. To prevent zones from becoming excessively large, reflecting the design of micro-transit systems as a complement rather than a replacement to public transit, we set a single-zone budget of $B_0 = 2$. Under this configuration, optimal solutions are expected to contain three to five zones, each less the size of one-fifth of a city. While the cost structure and coefficients are artificial, we are collaborating with our partner transit agency, CARTA, to calibrate these parameters in an ongoing micro-transit pilot deployment.

Experiments are run individually on a c6a.xlarge AWS instance with 4 CPUs and 8 GiB memory. For exact pricing, the ILP (\ref{ilp:pricing}) is solved with a 10-hour limit for CG and 1 hour per pricing problem. For the pricing heuristic, CG is limited to 20 minutes and each 2 minute per pricing problem. The number of random runs $R$ in Algorithm \ref{alg:pricing_heuristic} is set to $10$. All mathematical programs involved are solved by \href{https://www.gurobi.com/}{Gurobi}.

\subsection{Main Results}
The baseline algorithm is from \cite{hu2025optimal}, hereafter referred to as \textsc{CliqueGen}. Because the original implementation does not incorporate a global budget $B$ and requires specifying the number of zones $m$, we modify it to compute the operational cost of any candidate zone, defined as $\alpha \cdot D^2 + \beta$, alongside the clique construction process.

We first compare the computational times between \textsc{CliqueGen} and \textsc{CG}. From Table \ref{tab:cliquegen_cg_time}, we can see that \textsc{CliqueGen} only completes successfully for Miami and Boston with resolution-7 H3 hexagons. For all larger instances, \textsc{CliqueGen} terminates with the out-of-memory error in less than $1$ hour on the dedicated AWS instance. In contrast, \textsc{CG} + \textsc{ExactPricing} successfully completed all runs within the $10$ hours time limit. When the problem size scales from Miami to Boston, the computational time for the CG framework increases by only 3x (from 0.16 to 0.49 seconds), while \textsc{CliqueGen} increases by more than 100x (from 0.95 to 106 seconds). This dramatic difference demonstrates that \textsc{CliqueGen} suffers from severe scalability issues and is not effective for large cities.
\begin{table}[htbp!]
\centering
\caption{Computational times in seconds between \textsc{CliqueGen} and \textsc{CG} + \textsc{ExactPricing} across all \{city, H3 resolution\} configurations. OOM means the out-of-memory error.}
\label{tab:cliquegen_cg_time}
\begin{tabular}{>{\centering\arraybackslash}p{2.5cm} >{\centering\arraybackslash}p{1.5cm} >{\centering\arraybackslash}p{3cm}}
\toprule
\textbf{Instances} & \textbf{\textsc{CliqueGen}} & \textbf{\textsc{CG}} + \textbf{\textsc{Exact}} \textbf{\textsc{Pricing}} \\
\midrule
Miami-7   & 0.95  & 0.16 \\
Boston-7  & 106   & 0.49 \\
Atlanta-7 & OOM & 7.56 \\
Chattanooga-7   & OOM & 9.94 \\
Others          & OOM & Completed \\
\bottomrule
\end{tabular}
\end{table}

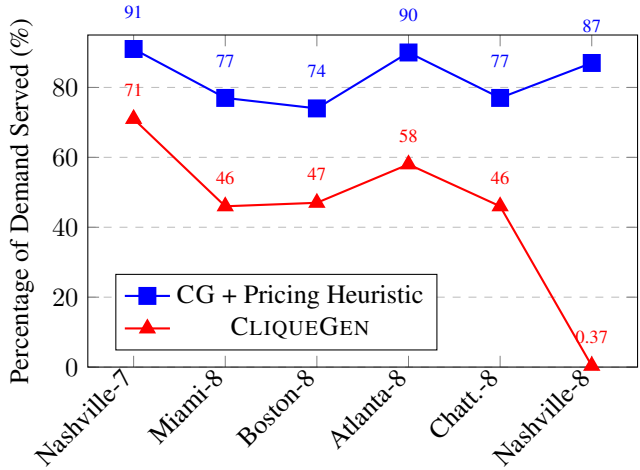
\begin{figure}[htbp!]
\centering
\begin{tikzpicture}
\begin{axis}[
    width=0.5\textwidth,
    height=6cm,
    ylabel={Percentage of Demand Served (\%)},
    ymin=0,
    ymax=95,
    xtick=data,
    xticklabels={Nashville-7, Miami-8, Boston-8, Atlanta-8, Chatt.-8, Nashville-8},
    x tick label style={rotate=45, anchor=east},
    legend style={at={(0.05,0.05)}, anchor=south west},
    ymajorgrids=true,
    grid style=dashed,
    mark size=3pt,
    every axis plot/.append style={thick}
]
\addplot[
    color=blue,
    mark=square*,
    nodes near coords,
    every node near coord/.append style={font=\scriptsize, anchor=south, yshift=8pt},
    point meta=explicit symbolic
    ]
    coordinates {
    (1,91)[91](2,77)[77](3,74)[74](4,90)[90](5,77)[77](6,87)[87]
    };
\addlegendentry{\textsc{CG} + Pricing Heuristic}
\addplot[
    color=red,
    mark=triangle*,
    nodes near coords,
    every node near coord/.append style={font=\scriptsize, anchor=north, yshift=17pt},
    point meta=explicit symbolic
    ]
    coordinates {
    (1,71)[71](2,46)[46](3,47)[47](4,58)[58](5,46)[46](6,0.37)[0.37]
    };
\addlegendentry{\textsc{CliqueGen}}
\end{axis}
\end{tikzpicture}
\caption{Comparison of demand coverage between \textsc{CliqueGen} and \textsc{CG} + \textsc{PricingHeuristic} within a 20-minute compute budget}
\label{fig:cliquegen_cg_comparison}
\end{figure}

Since both \textsc{CliqueGen} and \textsc{CG} + \textsc{PricingHeuristic} are anytime algorithms, we impose the same time limit of 1200 seconds on both methods and compare the best solutions obtained across all large instances with more than 200 hexagons. To ensure robustness, the results from \textsc{CG} are averaged over $5$ independent experiments. Table \ref{fig:cliquegen_cg_comparison} demonstrates that \textsc{CG} + \textsc{PricingHeuristic} significantly outperforms \textsc{CliqueGen}, achieving approximately 38\% higher demand coverage on average.


The performance gap is particularly pronounced in the largest instance: Nashville at resolution-8 with 1803 hexagons. Here, \textsc{CG} + \textsc{PricingHeuristic} leverages dual information to guide its search, exploring approximately 50 to 60 high-quality candidate zones and achieving 87\% demand coverage. In stark contrast, \textsc{CliqueGen} generates over 227,106 very small zones due to its iterative construction of cliques from smaller cliques, resulting in highly unstable outcome with only 0.37\% demand coverage.

This dramatic distinction highlights a fundamental advantage of CG's dual driven strategy: it uses economic signals from dual prices to prioritize zones that maximize marginal contribution to demand coverage.

\begin{figure}[htbp!]
    \centering
    \includegraphics[width=\linewidth]{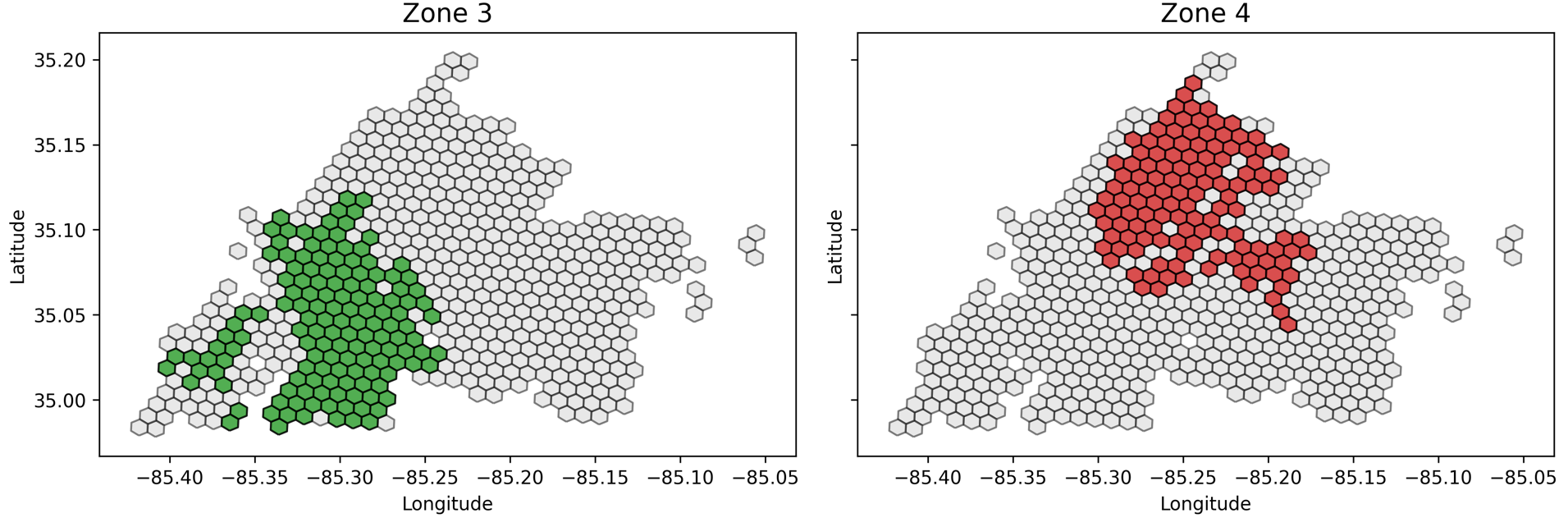}
    \caption{The micro-transit zones selected by \textsc{CG} + \textsc{PricingHeursitic} in Chattanooga at H3 resolution 8}
    \label{fig:optimal_zone}
\end{figure}

As an illustration, Figure~\ref{fig:optimal_zone} visualizes the zones selected by \textsc{CG} + \textsc{PricingHeuristic} for Chattanooga at resolution~8, which largely cover the high-demand areas shown in Figure~\ref{fig:out_demand_chatt}.
\begin{figure}[htbp!]
    \centering
    \includegraphics[width=\linewidth]{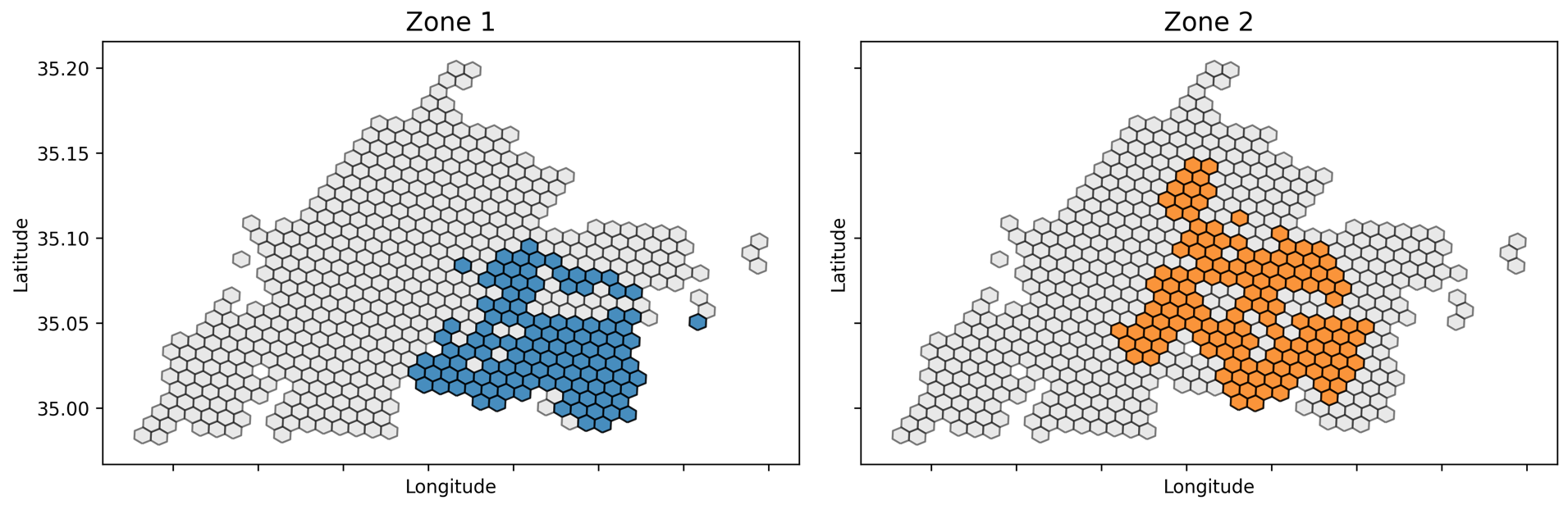}
\end{figure}
\vspace{-0.5cm}
\subsection{Exact vs. Heuristic Pricing}
We also compare the exact ILP formulation (ILP \ref{ilp:pricing}) and the heuristic (Algorithm \ref{alg:pricing_heuristic}) for solving the pricing problem. As shown in Table \ref{tab:exact_heuristic_pricing}, both approaches achieve nearly identical solution quality across the selected large instances. Except for Nashville at resolution 8 where the the exact ILP timeout without a result. They perform the same in all other instances. Given that the pricing heuristic runs approximately 30x faster than the exact method (from their time limits in Section \ref{sec:exp_setup}), it emerges as the preferred option for practical implementation.
\begin{table}[htbp!]
\centering
\caption{Comparison of demand coverage between Exact ILP formulation and pricing heuristic across five cities}
\label{tab:exact_heuristic_pricing}
\begin{tabular}{>{\centering\arraybackslash}p{2.5cm} >{\centering\arraybackslash}p{2.5cm} >{\centering\arraybackslash}p{2.5cm}}
\toprule
\textbf{Instances} & \textbf{Exact ILP} & \textbf{Pricing Heuristic} \\
\midrule
Boston-8      & 74.32\% & 74.21\% \\
Chattanooga-8 & 76.80\% & 76.70\% \\
Miami-8       & 79.92\% & 76.87\% \\
Atlanta-8     & 88.58\% & 89.76\% \\
Nashville-7   & 91.34\% & 91.32\% \\
Nashville-8   & Timeout & 90.97\% \\
\bottomrule
\end{tabular}
\end{table}

\vspace{-0.3cm}
\subsection{ILP vs. IQP in Exact Pricing}
We also compare the average computational times between the ILP formulation (\ref{ilp:pricing}) and the IQP formulation (\ref{qp:pricing}) for solving the pricing problem. As shown in Table \ref{tab:ilp_qp}, the IQP formulation consistently outperforms the ILP by 2x to 4x across all instances, despite quadratic programs being generally harder to solve. This gap likely stems from the excessive number of pairwise variables $y_{ij}$ introduced to the ILP, making the model less compact.
\begin{table}[htbp!]
\centering
\begin{tabular}{lccccc}
\toprule
 & \rotatebox{45}{Miami-8} & \rotatebox{45}{Boston-8} & \rotatebox{45}{Atlanta-8} & \rotatebox{45}{Chatt.-8} & \rotatebox{45}{Nash.-8} \\
\midrule
\textbf{ILP} & 350 & 486 & 1236 & $>$3600 & $>$3600 \\
\textbf{IQP} & 95 & 185 & 532 & 2016 & $>$3600 \\
\bottomrule
\end{tabular}
\caption{Average computation time in seconds for exact ILP and IQP pricing on selected large instances. All observations are capped by the 3600 seconds time limit.}
\label{tab:ilp_qp}
\end{table}
\vspace{-0.3cm}
\subsection{Experimental Details}
Due to space constraints, full experimental details are provided in Appendix D.

\section{Conclusion and Future Direction}
This paper proposes a Column Generation-based framework for the Micro-Transit Zoning Problem and demonstrates its effectiveness through extensive experiments on multiple large U.S. cities. Our solution approach, informed by practical insights from local city planners and transit agencies, such as our partner CARTA, has the potential to generate significant social impact by providing more sustainable and equitable micro-transit services.

Future work will consider a more general zoning model that jointly accounts for intra-zone and inter-zone demand, with cross-zone travel served by fixed-route transit, enabling micro-transit to operate as an effective first-and-last-mile connector to mass public transportation.


\newpage

\bibliographystyle{named}
\bibliography{ijcai26}

@article{lens2022zoning,
  title={Zoning, land use, and the reproduction of urban inequality},
  author={Lens, Michael C},
  journal={Annual Review of Sociology},
  volume={48},
  pages={421--439},
  year={2022},
  publisher={Annual Reviews}
}

@article{domingo2021effect,
  title={Effect of zoning plans on urban land-use change: A multi-scenario simulation for supporting sustainable urban growth},
  author={Domingo, Dar{\'\i}o and Palka, Ga{\"e}tan and Hersperger, Anna M},
  journal={Sustainable Cities and Society},
  volume={69},
  pages={102833},
  year={2021},
  publisher={Elsevier}
}

@article{dingil2025fostering,
  title={Fostering inclusive urban transportation in planning and policy-making: An umbrella review using ALARM methodology},
  author={Dingil, Ali Enes},
  journal={Sustainable Futures},
  volume={9},
  pages={100420},
  year={2025},
  publisher={Elsevier}
}

@article{chandra2020designing,
  title={Designing zoning systems for freight transportation planning: a GIS-based approach for automated zone design using public data sources},
  author={Chandra, Aithichya and Pani, Agnivesh and Sahu, Prasanta K},
  journal={Transportation Research Procedia},
  volume={48},
  pages={605--619},
  year={2020},
  publisher={Elsevier}
}

@article{chandra2021multi,
  title={A multi-objective genetic algorithm approach to design optimal zoning systems for freight transportation planning},
  author={Chandra, Aitichya and Sharath, MN and Pani, Agnivesh and Sahu, Prasanta K},
  journal={Journal of Transport Geography},
  volume={92},
  pages={103037},
  year={2021},
  publisher={Elsevier}
}

@article{liu2021mobility,
  title={Mobility service design via joint optimization of transit networks and demand-responsive services},
  author={Liu, Yining and Ouyang, Yanfeng},
  journal={Transportation Research Part B: Methodological},
  volume={151},
  pages={22--41},
  year={2021},
  publisher={Elsevier}
}

@article{bonner2023achieving,
  title={Achieving equitable outcomes through optimal design in the development of microtransit zones},
  author={Bonner, Taylor and Miller-Hooks, Elise},
  journal={Journal of Transport Geography},
  volume={112},
  pages={103696},
  year={2023},
  publisher={Elsevier}
}

@article{dantzig1960decomposition,
  title={Decomposition principle for linear programs},
  author={Dantzig, George B and Wolfe, Philip},
  journal={Operations research},
  volume={8},
  number={1},
  pages={101--111},
  year={1960},
  publisher={INFORMS}
}

@article{barnhart1998branch,
  title={Branch-and-price: Column generation for solving huge integer programs},
  author={Barnhart, Cynthia and Johnson, Ellis L and Nemhauser, George L and Savelsbergh, Martin WP and Vance, Pamela H},
  journal={Operations research},
  volume={46},
  number={3},
  pages={316--329},
  year={1998},
  publisher={INFORMS}
}

@article{gilmore1961linear,
  title={A linear programming approach to the cutting-stock problem},
  author={Gilmore, Paul C and Gomory, Ralph E},
  journal={Operations research},
  volume={9},
  number={6},
  pages={849--859},
  year={1961},
  publisher={INFORMS}
}

@book{desaulniers2006column,
  title={Column generation},
  author={Desaulniers, Guy and Desrosiers, Jacques and Solomon, Marius M},
  volume={5},
  year={2006},
  publisher={Springer Science \& Business Media}
}

@article{desrochers1992new,
  title={A new optimization algorithm for the vehicle routing problem with time windows},
  author={Desrochers, Martin and Desrosiers, Jacques and Solomon, Marius},
  journal={Operations research},
  volume={40},
  number={2},
  pages={342--354},
  year={1992},
  publisher={Informs}
}

@article{vance1997airline,
  title={Airline crew scheduling: A new formulation and decomposition algorithm},
  author={Vance, Pamela H and Barnhart, Cynthia and Johnson, Ellis L and Nemhauser, George L},
  journal={Operations Research},
  volume={45},
  number={2},
  pages={188--200},
  year={1997},
  publisher={INFORMS}
}

@article{van1999parallel,
  title={Parallel machine scheduling by column generation},
  author={Van Den Akker, Janna M and Hoogeveen, Jan A and van de Velde, Steef L},
  journal={Operations research},
  volume={47},
  number={6},
  pages={862--872},
  year={1999},
  publisher={INFORMS}
}

@inproceedings{riley2019column,
  title={Column generation for real-time ride-sharing operations},
  author={Riley, Connor and Legrain, Antoine and Van Hentenryck, Pascal},
  booktitle={International Conference on Integration of Constraint Programming, Artificial Intelligence, and Operations Research},
  pages={472--487},
  year={2019},
  organization={Springer}
}

@book{bertsimas1997introduction,
  title={Introduction to linear optimization},
  author={Bertsimas, Dimitris and Tsitsiklis, John N},
  volume={6},
  year={1997},
  publisher={Athena scientific Belmont, MA}
}

@article{hu2025optimal,
  title={Optimal Micro-Transit Zoning via Clique Generation and Integer Programming},
  author={Hu, Hins and Goswami, Rhea and Jiang, Hongyi and Samaranayake, Samitha},
  journal={arXiv preprint arXiv:2509.11445},
  year={2025}
}

@article{du1999stabilized,
  title={Stabilized column generation},
  author={Du Merle, Olivier and Villeneuve, Daniel and Desrosiers, Jacques and Hansen, Pierre},
  journal={Discrete Mathematics},
  volume={194},
  number={1-3},
  pages={229--237},
  year={1999},
  publisher={Elsevier}
}

@article{bertsimas2019optimizing,
  title={Optimizing schools’ start time and bus routes},
  author={Bertsimas, Dimitris and Delarue, Arthur and Martin, Sebastien},
  journal={Proceedings of the National Academy of Sciences},
  volume={116},
  number={13},
  pages={5943--5948},
  year={2019},
  publisher={National Academy of Sciences}
}

@article{sen2025moveod,
  title={MoveOD: Synthesizing Origin-Destination Commute Distribution from US Census Data},
  author={Sen, Rishav and Dubey, Abhishek and Mukhopadhyay, Ayan and Samaranayake, Samitha and Laszka, Aron},
  journal={arXiv preprint arXiv:2510.18858},
  year={2025}
}

@misc{gurobi,
  author = {{Gurobi Optimization, LLC}},
  title = {{Gurobi Optimizer Reference Manual}},
  year = 2024,
  url = "https://www.gurobi.com"
}

@book{hodges2010public,
  title={Public transportation's role in responding to climate change},
  author={Hodges, Tina},
  year={2010},
  publisher={Diane Publishing}
}

@article{miller2016public,
  title={Public transportation and sustainability: A review},
  author={Miller, Patrick and de Barros, Alexandre G and Kattan, Lina and Wirasinghe, SC},
  journal={KSCE Journal of Civil Engineering},
  volume={20},
  number={3},
  pages={1076--1083},
  year={2016},
  publisher={Springer}
}

@inproceedings{bardaka2024empathy,
  title={Empathy and AI: Achieving equitable microtransit for underserved communities},
  author={Bardaka, Eleni and Van\_Hentenryck, Pascal and Lee, Crystal Chen and Mayhorn, Christopher B and Monast, Kai and Samaranayake, Samitha and Singh, Munindar P},
  year={2024},
  organization={International Joint Conferences on Artificial Intelligence Organization}
}

@inproceedings{kooti2017analyzing,
  title={Analyzing Uber's ride-sharing economy},
  author={Kooti, Farshad and Grbovic, Mihajlo and Aiello, Luca Maria and Djuric, Nemanja and Radosavljevic, Vladan and Lerman, Kristina},
  booktitle={Proceedings of the 26th international conference on world wide web companion},
  pages={574--582},
  year={2017}
}

\newpage
\onecolumn
\appendix

\section{Visualization of Demand Pattern across Cities}
\begin{figure}[htbp!]
    \centering
    \includegraphics[width=\linewidth]{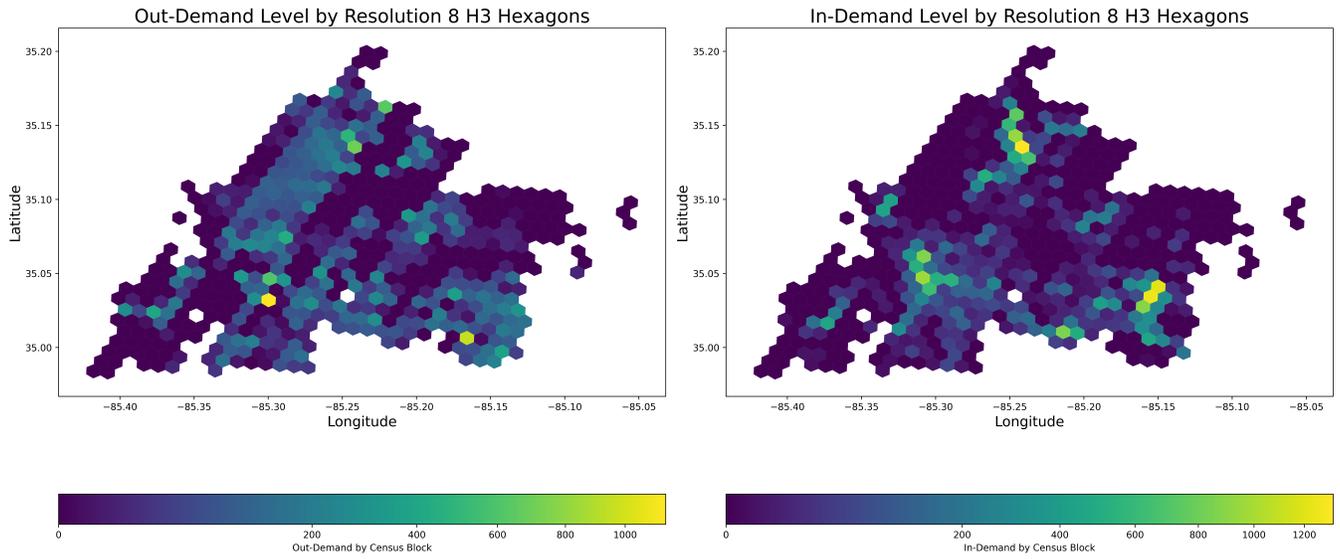}
    \caption{The number of trips in Chattanooga aggregated in resolution-8 H3 hexagons}
    \label{fig:chatt}
\end{figure}

\begin{figure}[htbp!]
    \centering
    \includegraphics[width=\linewidth]{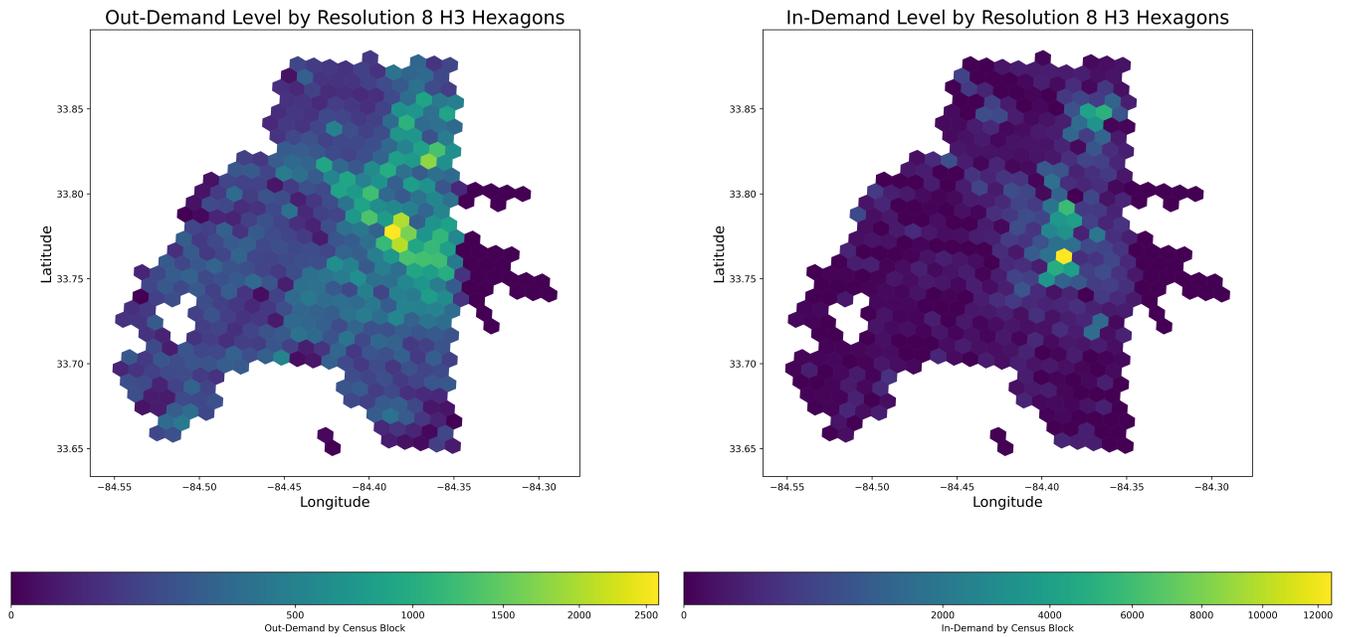}
    \caption{The number of trips in Atlanta aggregated in resolution-8 H3 hexagons}
    \label{fig:atlanta}
\end{figure}

\begin{figure}[htbp!]
    \centering
    \includegraphics[width=\linewidth]{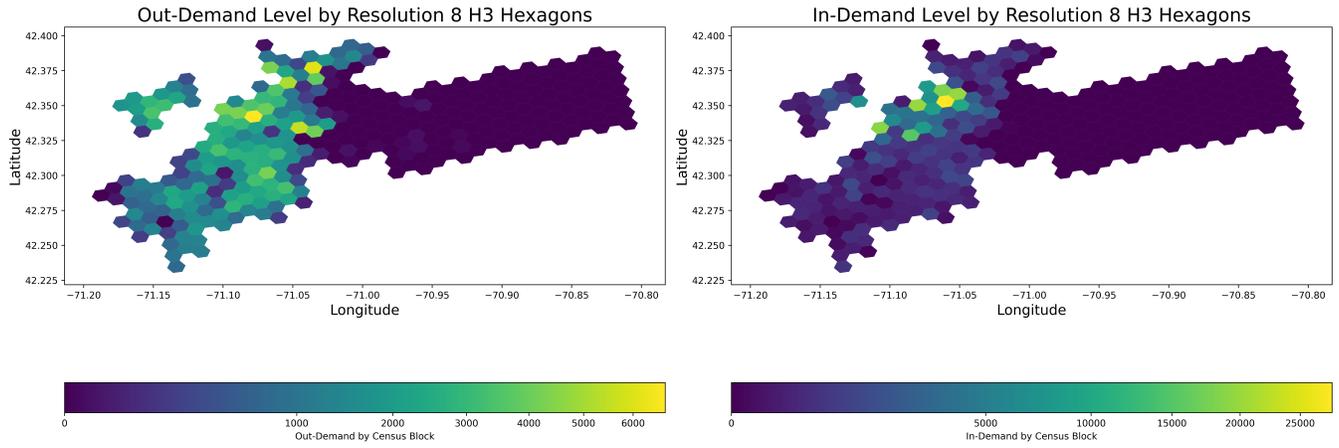}
    \caption{The number of trips in Boston aggregated in resolution-8 H3 hexagons}
    \label{fig:boston}
\end{figure}

\begin{figure}[htbp!]
    \centering
    \includegraphics[width=\linewidth]{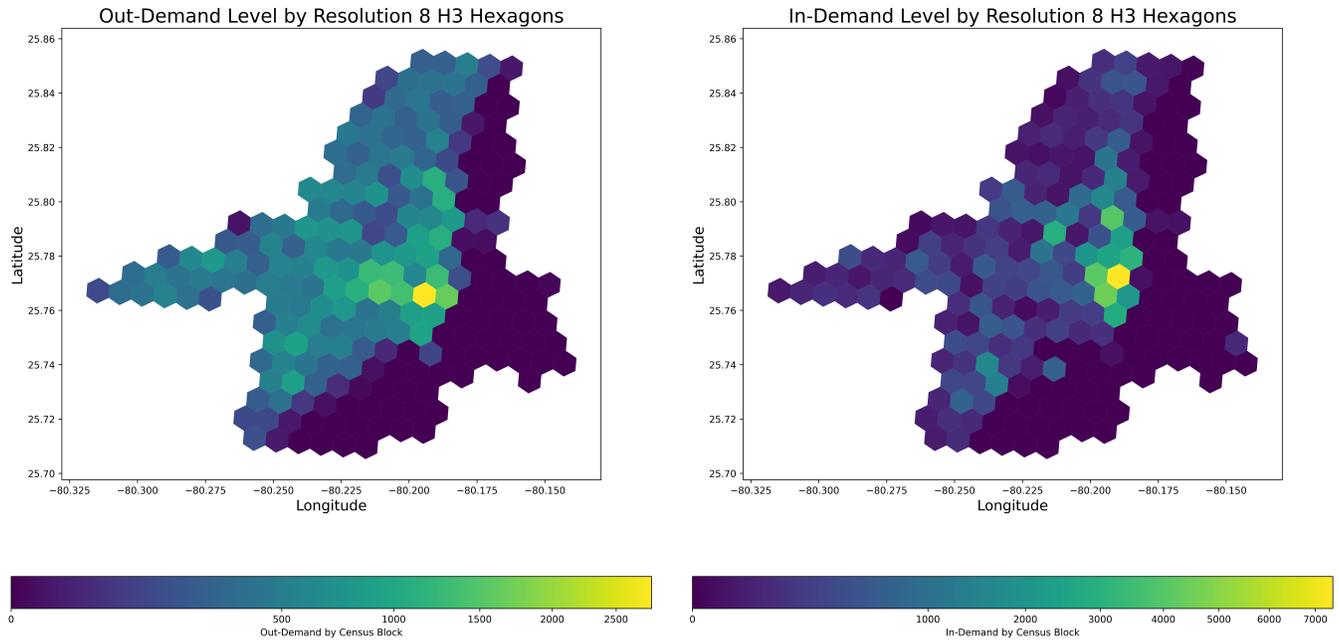}
    \caption{The number of trips in Miami aggregated in resolution-8 H3 hexagons}
    \label{fig:miami}
\end{figure}

\begin{figure}[htbp!]
    \centering
    \includegraphics[width=\linewidth]{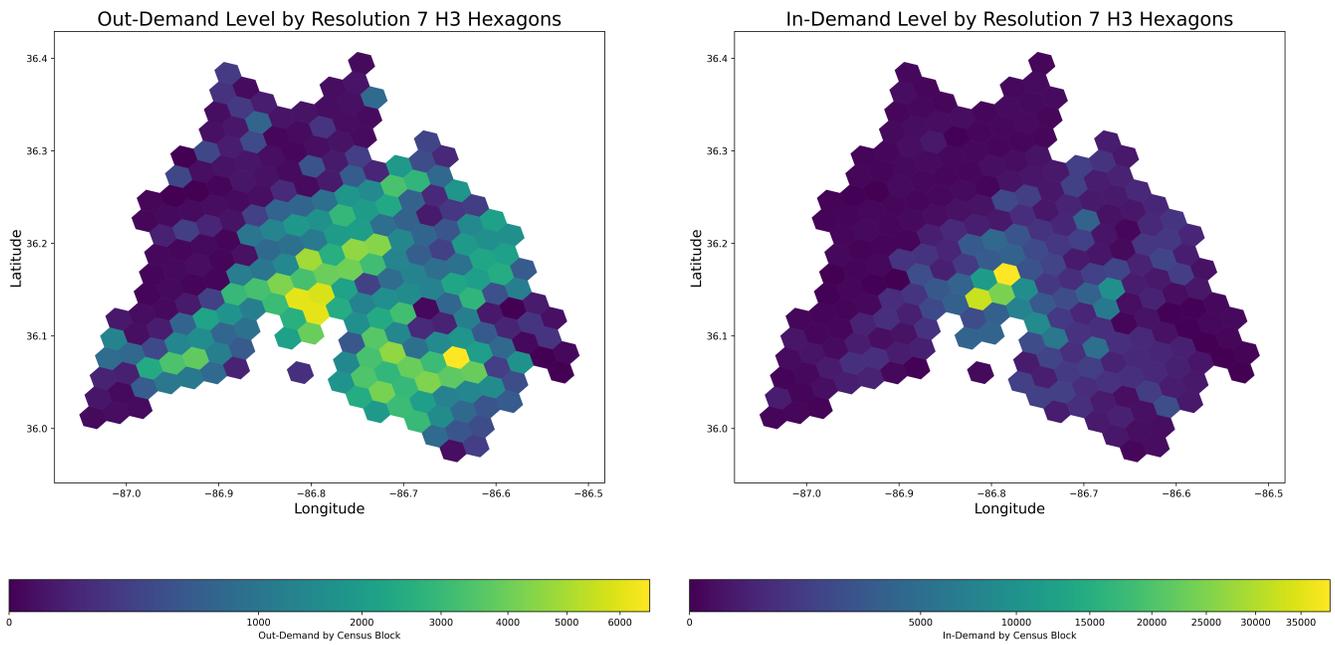}
    \caption{The number of trips in Nashville aggregated in resolution-7 H3 hexagons}
    \label{fig:nashville}
\end{figure}

Figures~\ref{fig:chatt}--\ref{fig:nashville} visualize the aggregated demand patterns for all cities discussed in the main text.


\newpage
\section{Diagrams for Root-Only Column Generation and Branch-n-Price}
\begin{figure}[htbp!]
    \centering
    \includegraphics[width=\linewidth]{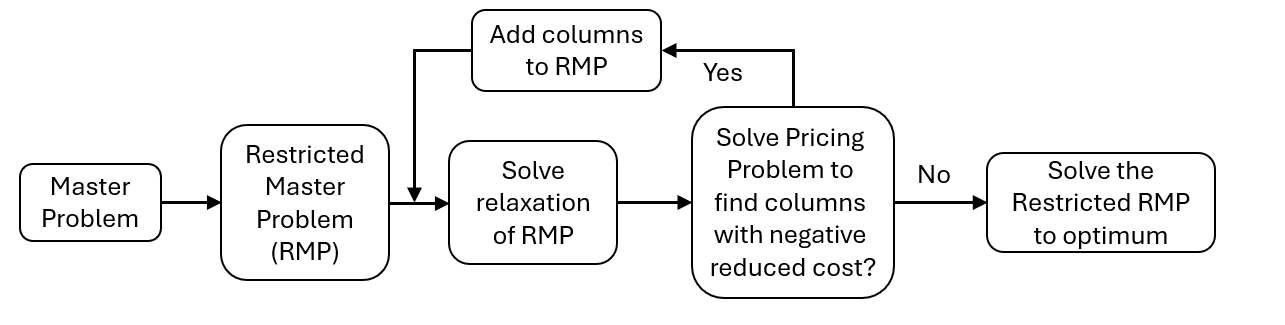}
    \caption{Diagram for CG on the LP relaxation to the restricted master problem}
    \label{fig:cg}
\end{figure}

\begin{figure}[htbp!]
    \centering
    \includegraphics[width=\linewidth]{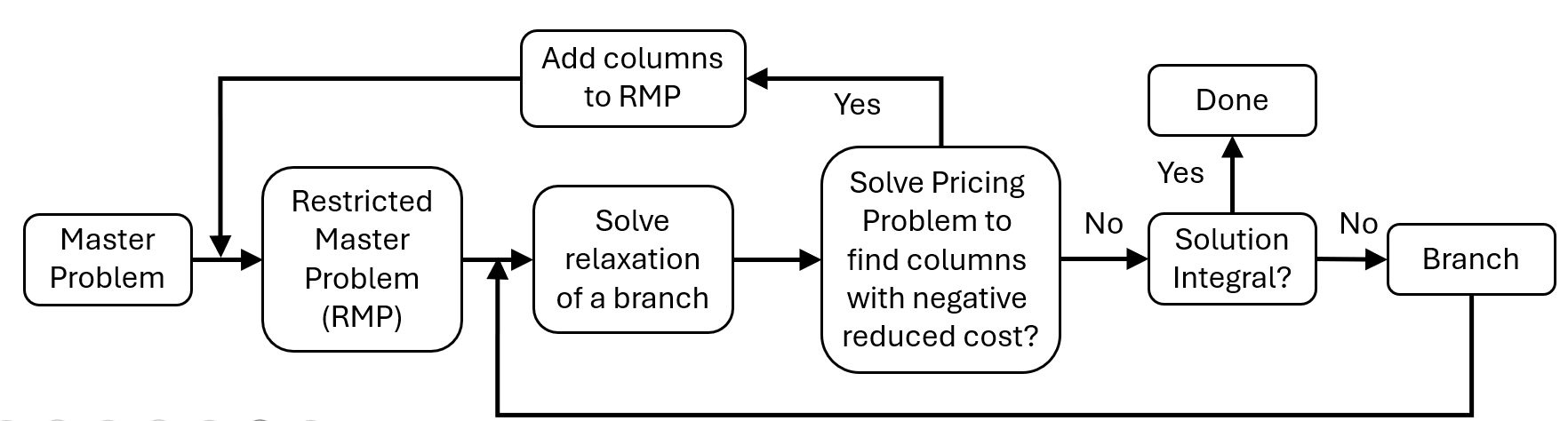}
    \caption{Diagram for branch-n-price}
    \label{fig:branch-n-price}
\end{figure}

Figures \ref{fig:cg} and \ref{fig:branch-n-price} illustrate root-only column generation and full Branch-and-Price, respectively.

\section{Demand Table}
Because city boundaries are fixed and our data pre-processor determines hexagon membership based on its centroid, the geo-fenced region for a city varies slightly across different H3 resolutions. Consequently, total travel demand differs between resolutions. Finer-resolution hexagons provide a more granular representation of the true spatial distribution of demand, and thus are preferable when computational resources allow.

\begin{table}[htbp!]
\centering
\caption{Number O-D trips aggregated under different H3 resolutions across multiple cities}
\label{tab:demand}
\begin{tabular}{ccc}
\toprule
\textbf{City} & \multicolumn{2}{c}{\textbf{Demand}} \\
\cmidrule(lr){2-3}
 & \textbf{Resolution 7} & \textbf{Resolution 8} \\
\midrule
Chattanooga  &  30480 &  31779 \\
Miami        &  67828 &  73013 \\
Atlanta      & 120316 & 124421 \\
Boston       & 207140 & 273272 \\
Nashville    & 290014 & 291139 \\
\bottomrule
\end{tabular}
\end{table}

\section{Total Demand Served by Different Methods across Cities}
CG with exact ILP pricing is run with a total time limit of 10 hours and a 1-hour limit per pricing problem. Both CliqueGen and CG with pricing heuristics are run under a 20-minute time limit, with a 2-minute limit per pricing problem applied to CG only.

\begin{table}[htbp!]
\centering
\caption{Total demand served (number of trips) across cities under two resolution levels.}
\label{tab:comprehensive_comparison}
\begin{tabular}{lccc}
\toprule
\textbf{City} & \textbf{CliqueGen} & \textbf{Exact ILP} & \textbf{Pricing Heuristic} \\
\midrule
\multicolumn{4}{c}{\textbf{Resolution 7}} \\
\midrule
Miami         & 54178       & 54178   & 54178   \\
Boston        & 154656      & 154656  & 154656  \\
Atlanta       & Not run yet & 105211  & 105211  \\
Chattanooga   & Not run yet & 22178   & 22178   \\
Nashville     & 204627      & 264908  & 264853  \\
\midrule
\multicolumn{4}{c}{\textbf{Resolution 8}} \\
\midrule
Miami         & 33878  & 58353   & 56125  \\
Boston        & 128258 & 203103  & 202790 \\
Atlanta       & 71715  & 110218  & 111674 \\
Chattanooga   & 14535  & 24407   & 24375  \\
Nashville     & 1084   & Timeout & 253306 \\
\bottomrule
\end{tabular}
\end{table}

\newpage
\section{Number of runs in pricing heuristic}
We observe that, under a fixed time limit, increasing the number of runs $R$ in the pricing heuristic (Algorithm 1) improves solution quality.

\begin{table}[htbp!]
\centering
\caption{Comparison of pricing heuristic performance with 5 runs vs 10 runs}
\label{tab:pricing_heuristic_runs}
\begin{tabular}{lcc}
\toprule
\textbf{City} & \textbf{5 Runs} & \textbf{10 Runs} \\
\midrule
\multicolumn{3}{c}{\textbf{Resolution 7}} \\
\midrule
Miami         & 54178   & 54178            \\
Boston        & 154656  & 154656           \\
Atlanta       & 105211  & 105211           \\
Chattanooga   & 22178   & 22178            \\
Nashville     & 224623  & \textbf{264853}  \\
\midrule
\multicolumn{3}{c}{\textbf{Resolution 8}} \\
\midrule
Miami         & 55906   & \textbf{56125}   \\
Boston        & 188343  & \textbf{202790}  \\
Atlanta       & 111343  & \textbf{111674}  \\
Chattanooga   & 24336   & \textbf{24375}   \\
Nashville     & 215644  & \textbf{253306}  \\
\bottomrule
\end{tabular}
\end{table}

\end{document}